\documentclass[reqno,12pt]{amsart}
\theoremstyle{plain}    
\newtheorem{thm}{Theorem}[section]
\theoremstyle{plain}    
\newtheorem{prop}[thm]{Proposition}
\theoremstyle{remark}

\newtheorem{example}[thm]{Example}
\newtheorem{examples}[thm]{Examples}

\usepackage{amscd,amssymb,comment,epic,euscript,graphics}

\newcommand\clspan{{\overline{\mathrm{span}}\,}}

\newcommand\Cpx{{\mathbf C}}

\newcommand\Gc{{\mathcal{G}}}

\newcommand\Ints{{\mathbf Z}}

\newcommand\Pt{{\widetilde P}}

\newcommand\qt{{\tilde q}}

\newcommand\Tcirc{{\mathbf T}}

\begin{document}

\pagestyle{myheadings}

 \title{Symmetric random walks on certain amalgamated free product groups}

 \author{Ken Dykema}

 \address{\hskip-\parindent
 Mathematisches Institut \\
 Westf\"alische Wilhelms--Universit\"at M\"unster \\
 Einsteinstr.\ 62 \\
 48149 M\"unster \\
 Germany}

 \address{\hskip-\parindent
 {\bf Permanent Address:} Department of Mathematics \\
 Texas A\&M University \\
 College Station TX 77843--3368, USA}
 \email{kdykema@math.tamu.edu}

 \thanks{Supported in part by a grant from the NSF and by the Alexander von Humboldt Foundation.
         The author also wishes to thank the Mathematics Institute of the
         Westf\"alische Wilhelms--Universit\"at
         M\"unster, for its generous hospitality during the author's year--long visit, when
         some of this research was conducted.}

 \date{4 October, 2004}

\begin{abstract}
We consider nearest--neighbor random walks on free products of finitely many copies of the integers
with amalgamation over nontrivial subgroups.
When all the subgroups have index two, we find the Green function of the random walks in terms of
complete elliptic integrals.
Our technique is to apply Voiculescu's operator--valued R--transform.
\end{abstract}

 \maketitle

 \markboth{\sc Random walks on amalgamated free products}{\sc Random walks on amalgamated free products}

\section{Introduction and description of results}

Symmetric random walks on groups have been much studied
since Kesten's classic paper~\cite{K} on the subject.
See, for example,~\cite{W:surv} and~\cite{W:book} and references therein.

Given a group $G$ with a finite, symmetric generating set $S$, we consider the
random walk on the associated Cayley graph of the group which starts
at the identity element $e$ and at each step
moves to any nearest neighbor with equal probability.
Let $p_n$ be the probability of return to $e$ on the $n$th step.
The {\em Green function} of the random walk is
\begin{equation}\label{eq:Green}
\Gc(z)=\sum_{n=0}^\infty p_n z^n
\end{equation}
and the {\em spectral radius} of the random walk is the reciprocal of the radius
of convergence of this power seris, i.e.
\[
r=\limsup_{n\to\infty}p_n^{1/n}.
\]

We will consider groups of the form $G=\Ints*_H\Ints$ or $G=(*_H)_1^N\Ints$,
which are amalgamated free products
of copies of the integers over a subgroup $H\cong\Ints$, where $H$ is embedded in the $j$th copy
of $\Ints$ as the subgroup of index $m_j$.
We will denote this group by $G=G_{m_1,m_2,\ldots,m_N}$.
Our generating sets will be $S=S_{m_1,\ldots,m_N}=\{a_1,a_1^{-1},\ldots,a_N,a_N^{-1}\}$,
where $a_j$ is a generator
of the $j$th copy of $\Ints$.
For convenience, we will call the corresponding random walk the {\em standard random walk}
on the group $G_{m_1,\ldots,m_N}$.
When $m_j=2$ for all $j$, we will  write the Green function of the random
walk in terms of Legendre's complete elliptic integrals.
In other cases, the Green function is equal to the integral of an algebraic function, and arbitarily many
terms of its power series expansion can be easily found.

It is known that the spectral radius of the standard random walk on $G=G_{m_1,\ldots,m_N}$ is an algebraic number.
Indeed, $H$ is a normal (in fact, central) subgroup of $G$ that is amenable.
By~\cite[Cor. 2]{K} and~\cite{K:amen}, the spectral radius of
the standard random walk on $G$ equals the spectral radius of the resulting random walk on the quotient group
$G/H$.
But
\[
G/H\cong(Z/m_1Z)*(Z/m_2Z)*\cdots*(Z/m_NZ)
\]
with the generating set $S_{m_1,\ldots,m_N}$ mapping to a union of generating sets for the cyclic groups $Z/m_jZ$.
The algebraicity of the spectral radius of this random walk is well known and has been proved by
different authors;  see~\cite{W:fp}, \cite{CS}, \cite{AK}.
See~\cite{CS} and~\cite{PW} for some results about random walks on other amalgamated free product groups. 

Our techniques rely on Voiculescu's operator--valued free probability theory,
and in particular on the operator--valued R--transform;
this is reviewed below, but see~\cite{V}, \cite{VDN} or~\cite{Sp} for more.

The paper is organized as follows.
In~\S\ref{sec:prelims} we review the elements of operator--valued free probability theory
that we will need, principally Voiculescu's R--transform.
In~\S\ref{sec:Green} we will show that the $B$--valued  Cauchy transform of the standard random
walk on $G_{m_1,\ldots,m_N}$ is an algebraic function, where $B=C^*_r(H)$ is the
reduced group C$^*$--algebra of the group over which we amalgamate.
We will also find the Green function when $m_1=\cdots=m_N=2$.

\vskip1ex
\noindent
{\em Acknowledgment.}  The author thanks Franz Lehner for informing him
that the spectral radius of these random walks is known, and how.

\section{Preliminaries}
\label{sec:prelims}

Let $A$ be a unital C$^*$--algebra and let
$B$ be a unital C$^*$--subalgebra of $A$ having a conditional expectation
$E:A\to B$.
The pair $(A,E)$ forms what is called a {\em $B$--valued noncommutative probability space.}
Given $T\in A$, the {\em $B$--valued Cauchy transform} of $T$ is the function
\[
 C^{(B)}_T(b)=\sum_{n=0}^\infty E((bT)^nb)=E((1-bT)^{-1}b)
\]
from a neighborhood of $0$ in $B$ into $B$.

We now review the means, devised by Voiculescu~\cite{V},
of finding the {\em $B$--valued R--transform} $R^{(B)}_T$ of $T$ from the $B$--valued
Cauchy transform.
The function $ C^{(B)}_T$ has an inverse $K^{(B)}_T=( C^{(B)}_T)^{\langle-1\rangle}$
with respect to composition.
Moreover, (see~\cite[Prop.\ 2.3]{Aa}), $K^{(B)}_T$ maps some neighborhood of $0$ in $B$
bijectively onto a neighborhood of $0$ in $B$ and maps invertible elements in this neighborhood
to invertible elements.
Then the R--transform of $T$ is
\[
R^{(B)}_T(b)=K^{(B)}_T(b)^{-1}-b^{-1}.
\]
This is defined for $b$ invertible and of small norm.
However, given that $ C^{(B)}_T$ is a power--series like sum of multilinear maps,
both $K^{(B)}_T$ and $R^{(B)}_T$ are seen to have a similar structure;
see also the combinitorial description in~\cite{Sp}.
The above definition of $R^{(B)}_T$ on invertible elements actually determines it on all
of a neighborhood of $0$ in $B$.

Suppose for every $i$ in an index set $I$, $A_i\subseteq A$ is a C$^*$--subalgebra of $A$ with $B\subseteq A_i$.
The family $(A_i)_{i\in I}$ is said to be {\em free} if $E(a_1\ldots a_n)=0$ whenever $a_j\in A_{i_j}\cap\ker E$
and $i_1\ne i_2,\,i_2\ne i_3,\ldots,i_{n-1}\ne i_n$.
A family $(T_i)_{i\in I}$ of elements of $A$ is said to be free if the family
$(C^*(B\cup\{T_i\}))_{i\in I}$ of C$^*$--subalgebras is free.

\begin{thm}[Voiculescu~\cite{V}]\label{thm:Rtrans}
In a $B$--valued noncommutative probability space $(A,E)$, suppose $T_i\in A$ are such that $(T_i)_{i=1}^N$ is
a free family.
Then
\[
R^{(B)}_{T_1+\cdots+T_N}=R^{(B)}_{T_1}+\cdots+R^{(B)}_{T_N}.
\]
\end{thm}

We will now consider how these results from free proability theory may be applied to
the study of random walks on amalgamated free products of groups.
Given a group $G$, let $\lambda=\lambda^{(G)}$ be the left regular representation of $G$ as unitary operators
on $\ell^2(G)$, extended linearly to a $*$--representation of the complex group algebra $\Cpx[G]$.
The {\em reduced group C$^*$--algebra} is
\[
\clspan\{\lambda(g)\mid g\in G\}.
\]
(We will sometimes write $\lambda_g$ instead of $\lambda(g)$.)
Note that the canonical trace on $\Cpx[G]$, which extracts the coefficient of the identity element,
extends to the tracial state $\tau=\tau_G$ on $C^*_r(G)$ given by $\tau(x)=\langle x\delta_e,\delta_e\rangle$.
Thus, $(C^*_r(G),\tau)$ is a $\Cpx$--valued noncommutative probability space.
If $H$ is a subgroup of $G$, then $\overline{\lambda^{(G)}(\Cpx[H])}$ is isomorphic to $C^*_r(H)$, and will
be denoted as such.
The projection $\ell^2(G)\to\ell^2(H)$ impliments the {\em canonical conditional expecation}
$E:C^*_r(G)\to C^*_r(H)$, which satisfies
\[
E(\lambda(g))=
\begin{cases}
\lambda(g),&g\in H \\
0,&\text{otherwise}.  
\end{cases}
\]
Thus, $(C^*_r(G),E)$ is a $C^*_r(H)$--valued noncommutative probability space.
Note that $E$ is $\tau$--preserving.

Suppose $G=(*_H)_{i=1}^NG_i$ is a free product of groups $G_i$ with amalgamation over subgroups $H\subseteq G_i$.
For ease of writing, we will take $N=2$ and write $G=G_1*_H G_2$,
though similar considerations apply for general $N<\infty$.
The group inclusions $G_i\subseteq G$ give rise to C$^*$--subalgebras $C^*_r(G_i)\subseteq C^*_r(G)$,
and the pair $(C^*_r(G_1),C^*_r(G_2))$ is
free with respect to the canonical conditional expectation $E:C^*_r(G)\to C^*_r(H)$.
Suppose $S_i$ is a finite, symmetric generating set for $G_i$.
Then $S=S_1\cup S_2$ a finite, symmetric generating set for $G$.
Let $T_i=\sum_{a\in S_i}\lambda(a)\in C^*_r(G_i)\subseteq C^*_r(G)$ and let $T=T_1+T_2$.
The {\em adjacency operator} $T$ is related to the random walk associated to $S$ by
$\tau(T^n)=|S|^np_n$ and the $\Cpx$--valued Cauchy transform
\[
 C^{(\Cpx)}_T(\zeta)=\sum_{n=0}^\infty\tau(T^n)\zeta^{n+1}
\]
is related to the Green function $\Gc$ of this random walk~\eqref{eq:Green} by
\begin{equation}\label{eq:CGr}
C^{(\Cpx)}_T(\zeta)=\zeta\Gc(|S|\zeta).
\end{equation}
Let $B=C^*_r(H)$,
We have that $ C^{(\Cpx)}_T$ is $\tau\circ C^{(B)}_T$ restricted to scalar multiples of the identity operator.

What follows, then, is a strategy for finding the Green function of this random walk on $G_1*_HG_2$.
If the $B$--valued Cauchy transforms $ C^{(B)}_{T_1}$ and $ C^{(B)}_{T_2}$ are known, then the
$B$--valued R--transforms of $T_1$ and $T_2$ can be found and used, with Voiculescu's Theorem~\ref{thm:Rtrans}
to compute the $B$--valued R--transform of $T$, from with the $B$--valued Cauchy transform of $T$ can be found.
Composing with $\tau$ then yields the $\Cpx$--valued Cauchy transform of $T$.

\section{Green functions}
\label{sec:Green}

In this section, we consider standard random walks on the amalgamated
free product groups $G_{m_1,\ldots,m_N}$.
We find the Green function $\Gc_{m_1,\ldots,m_N}$ of this
random walk when $m_j=2$ for all $j$, and we show how to derive information about the Green function
in other cases.

\begin{prop}\label{prop:ZGr}
Let $n\in\{2,3,4,\ldots\}$ and consider the subgroup
$n\Ints=H\subseteq G=\Ints$.
Let $E$ be the canonical conditional expectation from $A=C^*_r(G)$ onto the
subalgebra $B=C^*_r(H)$.
Consider the adjacency operator $T=\lambda_1+\lambda_{-1}\in A$.
Then the $B$--valued Cauchy transform of $T$ is
\begin{equation}\label{eq:ZGr}
 C^{(B)}_T(b)=\frac{bp(b)}{q(b)-b^n(\lambda_{n}+\lambda_{-n})},
\end{equation}
where $p$ and $q$ are polynomials with integer coefficients, each with
constant term equal to $1$ and with $\deg(p)\le n-1$ and $\deg(q)\le n$.
\end{prop}
\begin{proof}
Using the Fourier transform, $A=C^*_r(\Ints)$ is seen to be isomorphic to
$C(\Tcirc)$, the algebra of all continuous functions on the circle, and we
henceforth make this identification of $A$ with $C(\Tcirc)$.
Thus $\lambda_k\in A$ is identified with the function that is the map
$\Tcirc\ni z\mapsto z^k$.
The subalgebra $B=C^*_r(n\Ints)$, is identified with
the set of functions invariant under rotation of the domain $\Tcirc$ by angle $2\pi/n$,
and for every $f\in A$,
\[
(Ef)(z)=\frac1n\sum_{k=0}^{n-1}f(\omega^kz),
\]
where $\omega=\omega_n=\exp(2\pi i/n)$.
For $b\in B$ of sufficiently small norm,
\begin{align}
\notag
 C^{(B)}_T(b)&=E(b(1-Tb)^{-1})=
\frac bn\sum_{k=0}^{n-1}
\frac1{1-(\omega^k\lambda_1+\omega^{-k}\lambda_{-1})b} \\[4ex]
&=\bigg(\frac bn\sum_{k=0}^{n-1}
\prod_{\substack{j=0\\j\ne k}}^{n-1}
\big(1-(\omega^j\lambda_1+\omega^{-j}\lambda_{-1})b\big)\bigg)
\bigg/\bigg(\prod_{j=0}^{n-1}
\big(1-(\omega^j\lambda_1+\omega^{-j}\lambda_{-1})b\big)\bigg). \label{eq:GnNumDen}
\end{align}

Consider the ring
\[
R=\{\sum_{k\in\Ints}a_k\lambda_k\mid a_k\in\Ints[\omega],\text{ all but finitely many }a_k=0\},
\]
which is the group ring of $\Ints$ with coefficients from $\Ints[\omega]$.
Consider the denominator $Q$ of~\eqref{eq:GnNumDen} as a polynomial in variable $b$ with coefficients
from $R$.
Then $Q$ is of degree $\le n$ and has constant term equal to $1$.
The coefficient of $b^k$, $1\le k\le n$, is of the form
\begin{equation}\label{eq:bkcoeff}
a_{-k}\lambda_{-k}+a_{-k+1}\lambda_{-k+1}+\cdots+a_k\lambda_k
\end{equation}
with all $a_j\in\Ints[\omega]$.
We see from~\eqref{eq:GnNumDen}
that $Q$ is invariant under the automorphism of $R$ given by $\lambda_j\mapsto\omega^j\lambda_j$,
($j\in\Ints$).
So only $a_0$ can be nonzero in~\eqref{eq:bkcoeff} if $1\le k\le n-1$,
while with $k=n$, we find that the coefficient of $b^n$ in $Q$ is of the form
\begin{equation}\label{eq:bncoeff}
a_{-n}\lambda_{-n}+a_0+a_n\lambda_n.
\end{equation}
But since the coefficient of $b_n$ equals
\[
(-1)^n\prod_{j=0}^{n-1}(\omega^{-j}\lambda_{-1}+\omega_j\lambda_1),
\]
we see
\[
a_n=(-1)^n\prod_{j=0}^{n-1}\omega^j=(-1)^n\omega^{(n-1)n/2}=-1
\]
and $a_{-n}=a_n^{-1}=-1$.
Finally, $Q$ is invariant under the transformation $\omega\mapsto\omega^d$ whenever $d$
is relatively prime to $n$.
So by the fundamental theorem of Galois theory, we get $a_0\in\Ints$ in~\eqref{eq:bncoeff}
and the coefficient of $b^k$ is an integer for $1\le k\le n-1$.
Let $q$ be the polynomial given by $Q(b)=q(b)-b^n(\lambda_{-n}+\lambda_n)$.

By the same reasoning as for $Q$, we see that the numerator, $P$, of~\eqref{eq:GnNumDen} is
equal to $\frac1nb\Pt(b)$, where $\Pt$ is a polynomial of degree $\le n-1$ having
integer coefficients, and where the constant coefficient of $\Pt$ is equal to $n$.
We need only show that the coefficients of $\frac1n\Pt$ are all integers.
Equating two descriptions of $ C^{(B)}_T(b)$, we get
\[
\sum_{k=0}^\infty E(b(Tb)^k)=\frac{\frac bn\Pt(b)}{1-(\qt(b)+b^n(\lambda_{-n}+\lambda_n))}
=\frac bn\Pt(b)\sum_{k=0}^\infty(\qt(b)+b^n(\lambda_{-n}+\lambda_n))^k,
\]
where $\qt(b)=1-q(b)$.
Since $E(b(Tb)^k)$ can be written as a linear combination of $\{\lambda_j\mid j\in n\Ints,\,|j|\le k\}$
with coefficients from $\Ints[b]$ and since $\qt$ has integer coefficients and zero constant coefficient,
we conclude that $\frac 1n\Pt$ has integer coefficients.
\end{proof}

\begin{examples}\label{exs:Gn}
Let $ C_n$ denote the $B$--valued Cauchy transform considered
in equation~\eqref{eq:ZGr} of Proposition~\ref{prop:ZGr}.
Using the formula~\eqref{eq:GnNumDen}, we find
\begin{alignat*}{2}
 C_2(b)&=\frac b{1 - 2b^2 - b^2(\lambda_{-2} + \lambda_{2})} &\qquad
 C_3(b)&=\frac{b - b^3}{1 - 3b^2 - b^3(\lambda_{-3}+ \lambda_{3})} \\[1ex]
 C_4(b)&=\frac{b - 2b^3}{1 - 4b^2 + 2b^4 - b^4(\lambda_{-4}+\lambda_{4})} &
 C_5(b)&=\frac{b - 3b^3 + b^5}{1 - 5b^2 + 5b^4 - b^5(\lambda_{-5}+\lambda_{5})}
\end{alignat*}
\end{examples}

From Proposition~\ref{prop:ZGr} and the procedure for obtaining the R--transform
as described in~\S\ref{sec:prelims}, we see that the $C^*_r(H)$--valued R--transform of the operator $T$
considered above is an algebraic function.
A precise formulation is below.

\begin{prop}\label{prop:ZR}
Let $n\in\{2,3,4,\ldots\}$,
let $B=C^*_r(n\Ints)\subseteq A=C^*_r(\Ints)$, $E:A\to B$ and $T$ be as in Proposition~\ref{prop:ZGr}.
Let $R=R^{(B)}_T$ be the $B$--valued R--transform of $T$.
Then there is an irreducible polynomial $Q_n$ in three variables and with integer coefficients
such that
\[
Q_n(R(b),b,\xi)=0,
\]
where $\xi=\lambda_{-p}+\lambda_p$, for $p$ a generator of $H$.
\end{prop}

\begin{examples}\label{exs:Rn}
Below are listed some of the irreducible polynomials $Q_n=Q_n(R,b,\xi)$ from
Proposition~\ref{prop:ZR}.
\begin{align*}
Q_2&=bR^2+R-b(2+\xi) \\
Q_3&=b^2R^3+2bR^2+(1-3b^2)R-b(2+ b\xi) \\
Q_4&=b^3R^4+3b^2R^3+b(3-4b^2)R^2+(1-6b^2)R-b(2-2b^2+b^2\xi) \\
Q_5&=b^4R^5+4b^3R^4+(6b^2-5b^4)R^3+(4b-12b^3)R^2+(1-9b^2+5b^4)R \\
   &\quad-b(2-4b^2+b^3\xi) 
\end{align*}
\end{examples}

\begin{prop}\label{prop:Z*ZRGr}
Let $G=G_{m_1,\ldots,m_N}$
with generating set $S=S_{m_1,\ldots,m_N}$ be
as described in the introduction.
Let $B=C^*_r(H)\subseteq C^*_r(G)=A$ equipped with the canonical conditional expectation $E:A\to B$.
Let $T=\sum_{a\in S}\lambda_a\in A$ be the adjacency operator for the standard random walk.
Let $ C= C^{(B)}_T$ and $R=R^{(B)}_T$ be the $B$--valued Cauchy transform
and R--transform of $T$, respectively.
Then there is are irreducible polynomials $P_{m_1,\ldots,m_N}$ and $Q_{m_1,\ldots,m_N}$,
each in three variables and with integer coefficients, such that
\[
P_{m_1,\ldots,m_N}( C(b),b,\xi)=0,\qquad Q_{m_1,\ldots,m_N}(R(b),b,\xi)=0,
\]
where $\xi=\lambda_{-p}+\lambda_p$,
for $p$ a generator of $H$.
\end{prop}
\begin{proof}
We have $S=\{a_1^{-1},a_1,\ldots,a_N^{-1},a_N\}$ where $a_j$ is a generator of the $j$th
copy of $\Ints$ in the amalgamated free product $G=(*_H)_1^N\Ints$,
and $T=T_1+\cdots+T_N$, where $T_j=\lambda_{a_j^{-1}}+\lambda_{a_j}$.
By additivity of the $B$--valued R--transform (Theorem~\ref{thm:Rtrans}),
\[
R^{(B)}_T=R^{(B)}_{T_1}+\cdots+R^{(B)}_{T_N}.
\]
By Proposition~\ref{prop:ZR}, each $R^{(B)}_{T_j}(b)$ is an algebraic
function of $b$ and $\xi$, being the root of the polynomial with integer coefficients,
so the same is true for $R^{(B)}_T$.
Now the procedure for finding $ C^{(B)}_T$ from $R^{(B)}_T$ yields the polynomial $P_{m_1,\ldots,m_N}$
from $Q_{m_1,\ldots,m_N}$.
\end{proof}

The polynomials $Q_{m_1,\ldots,m_N}$ and $P_{m_1,\ldots,m_N}$ are easily found, as is illustrated
in the following three examples.
In the first two of these examples, we are able to write explicitly the Green functions of the random walks.

\begin{example}\label{ex:22}
Consider the case $G=G_{2,2}$.
Note that $G$ is an amenable group.
Then from $Q_2$ of Examples~\ref{exs:Rn}, we get immediately
\[
Q_{2,2}(R,b,\xi)=bR^2+2R-4b(2+\xi)
\]
and
\begin{equation}\label{eq:G22}
P_{2,2}( C,b,\xi)=(1 - 8b^2 - 4b^2\xi) C^2-b^2.
\end{equation}
Letting $ C^{(B)}_{2,2}$ denote the $B$--valued Cauchy transform of the adjacency operator $T$,
namely the quantity $ C$ in~\eqref{eq:G22} above, and
using the asymptotic behavior $ C^{(B)}_{2,2}(b)=b+O(\|b\|^2)$ as $\|b\|\to0$ to choose the
branch of the square root, we get
\[
 C^{(B)}_{2,2}(b)=\frac b{\sqrt{1-4b^2(2+\xi)}}.
\]
Letting $ C_{2,2}$ denote the scalar--valued Cauchy transform of $T$, we have for $\zeta\in\Cpx$,
\[
 C_{2,2}(\zeta)=\tau\circ C^{(B)}_{2,2}(\zeta),
\]
where $\tau$ is the canonical trace on $C^*_r(H)$.
Thus, taking $|\zeta|$ small,
\begin{align}
\notag
 C_{2,2}(\zeta)
&=\frac\zeta{2\pi i}\int_{|z|=1}
 \frac{dz}{z\sqrt{1-4\zeta^2(z^{-1}+z+2)}} \\
\notag
&=\frac\zeta{2\pi i}\int_{|z|=1}
 \frac{dz}{\sqrt{z(-4\zeta^2z^2+(1-8\zeta^2)z-4\zeta^2)}} \\
&=\frac1{4\pi i}\int_{|z|=1}
 \frac{dz}{\sqrt{z(z-z_1)}\,\sqrt{(z_2-z)}}, \label{eq:G22int}
\end{align}
where
\begin{align*}
z_1=z_1(\zeta)&=\frac{1-8\zeta^2-\sqrt{1-16\zeta^2}}{8\zeta^2}
 =4\zeta^2+O(|\zeta|^4) \\ \vspace{1ex}
z_2=z_2(\zeta)&=\frac{1-8\zeta^2+\sqrt{1-16\zeta^2}}{8\zeta^2}
 =\frac1{4\zeta^2}+O(1)
\end{align*}
with the indicated asymptotics as $\zeta\to0$
and where in~\eqref{eq:G22int}, the branch of $\sqrt{z_2-z}$ close to $1/2\zeta$ is chosen.
Replacing the contour $|z|=1$ by the contour drawn in Figure~\ref{fig:contour},
\begin{figure}[bht]
\begin{picture}(140,100)(-20,-20)
\put(0,0){\circle*{3}}
\put(100,60){\circle*{3}}
\put(0,0){\circle{40}}
\put(100,60){\circle{40}}
\put(17.15,10.29){\circle*{3}}
\drawline(17.15,10.29)(82.85,49.71)
\put(82.85,49.71){\circle*{3}}
\put(-3,-12){$0$}
\drawline(-17.15,-10.29)(-17.15,-2.29)
\drawline(-17.15,-10.29)(-24.21,-6.52)
\put(97,66){$z_1$}
\drawline(117.15,70.29)(117.15,62,29)
\drawline(117.15,70.29)(124.21,66,52)
\drawline(60,36)(56.24,28.94)
\drawline(40,24)(43.76,31.16)
\end{picture}
\caption{Countour used in evaluating the integral~\eqref{eq:G22int}}
\label{fig:contour}
\end{figure}
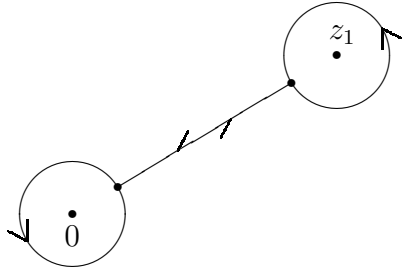
and letting the circles shrink to the point of disappearing,
we get
\begin{align}
\notag
 C_{2,2}(\zeta)&=\frac1{2\pi\sqrt{z_2}}\int_0^1\frac{dt}{\sqrt{t(1-t)}\sqrt{1-\frac{z_1}{z_2}t}}
=\frac1{\pi\sqrt{z_2}}\int_0^{\pi/2}\frac{d\phi}{\sqrt{1-\frac{z_1}{z_2}\sin^2\phi}} \\
&=\frac1{\pi\sqrt{z_2(\zeta)}}F_1(\sqrt{\tfrac{z_1(\zeta)}{z_2(\zeta)}}), \label{eq:F1}
\end{align}
where we have made the change of variables $t=\sin^2\phi$,
where $F_1$ is Legendre's complete elliptic integral of the first kind (see~\cite{Cayley})
and where we take the branch of $\sqrt{z_2(\zeta)}$ that is close to $1/(2\zeta)$;
it is not necessary to specify
the branch of $\sqrt{z_1(\zeta)/z_2(\zeta)}$.
Therefore, the Green function of the standard random walk on on $G_{2,2}$ is
\[
\Gc_{2,2}(z)=\frac4{\pi\sqrt{2-z^2+2\sqrt{1-z^2}}}
\,F_1\bigg(\sqrt{\frac{2-z^2-2\sqrt{1-z^2}}{2-z^2+2\sqrt{1-z^2}}}\bigg).
\]
\end{example}

\begin{example}\label{ex:2n2}
Let $G=G_{2,\ldots,2}=(*_\Ints)_1^N\Ints$ be the free product of $N\ge3$ copies of the integers
with amalgamation over their index--two subgroups.
As in the previous example and using analogous notation, we find
\begin{align}
\notag
Q_{2,\ldots,2}&=bR^2 + NR - bN^2(2 + \xi)   \\
P_{2,\ldots,2}&=(1 - b^2N^2(2+\xi)) C^2 + b(N-2) C - b^2(N-1). \label{eq:P2n2}
\end{align}
Solving for $ C$ gives
\[
 C^{(B)}_{2,\ldots,2}(b)=\bigg(\frac{2-N+N\sqrt{1-4(N-1)b^2(2+\xi)}}{2(1-N^2b^2(2+\xi))}\bigg)b,
\]
and integrating yields the $\Cpx$--valued Cauchy transform
\begin{align}
\notag
 C_{2,\ldots,2}(\zeta)&=\frac{\zeta}{2\pi i}\int_{|z|=1}
\frac{2-N+N\sqrt{1-4(N-1)\zeta^2(2+z^{-1}+z)}}{2z(1-N^2\zeta^2(2+z^{-1}+z))}dz \\[2ex]
&=\frac{(2-N)\zeta}{4\pi i}\int_{|z|=1}\frac{dz}{-N^2\zeta^2z^2+(1-2N^2\zeta^2)z-N^2\zeta^2} 
\label{eq:G2n2_1stint} \\[2ex]
&\quad+\frac{N\zeta}{4\pi i}\int_{|z|=1}
 \frac{\sqrt{1-4(N-1)\zeta^2(2+z^{-1}+z)}}{-N^2\zeta^2z^2+(1-2N^2\zeta^2)z-N^2\zeta^2}dz,
\label{eq:G2n2_2ndint}
\end{align}
for $\zeta$ sufficiently small.
The denominator in the integrals~\eqref{eq:G2n2_1stint} and~\eqref{eq:G2n2_2ndint} has roots
\begin{align}
z_3=z_3(\zeta)&=\frac{1-2N^2\zeta^2-\sqrt{1-4N^2\zeta^2}}{2N^2\zeta^2}
 =N^2\zeta^2+O(|\zeta|^4) \label{eq:z3} \\
z_4=z_4(\zeta)&=\frac{1-2N^2\zeta^2+\sqrt{1-4N^2\zeta^2}}{2N^2\zeta^2}
 =\frac1{N^2\zeta^2}+O(1). \label{eq:z4}
\end{align}
The value of the first term~\eqref{eq:G2n2_1stint} is, thus,
\begin{align}
\notag
\frac{(2-N)\zeta}{4\pi i}\int_{|z|=1}\frac{dz}{N^2\zeta^2(z-z_3)(z_4-z)}
&=\frac{(2-N)\zeta}{2}\;\frac1{N^2\zeta^2(z_4-z_3)} \\
&=\frac{(2-N)\zeta}{2\sqrt{1-4N^2\zeta^2}}. \label{eq:2-N}
\end{align}
The second term~\eqref{eq:G2n2_2ndint} equals
\begin{align}
\notag
\frac{N\zeta}{4\pi i}\int_{|z|=1}&\frac{\sqrt{\frac1z(z-4(N-1)\zeta^2(2z+1+z^2))}}
 {N^2\zeta^2(z-z_3)(z_4-z)}\,dz \\
&=\frac1{4N\zeta\pi i}\int_{|z|=1}\frac{p(z)}
 {(z-z_3)(z_4-z)\sqrt{z\,p(z)}}\, dz \label{eq:2ndterm}
\end{align}
where
\[
p(z)=z-4(N-1)\zeta^2(2z+1+z^2)=-4(N-1)\zeta^2z^2+(1-8(N-1)\zeta^2)z-4(N-1)\zeta^2.
\]
The roots of $p(z)$ are
\begin{align}
z_5=z_5(\zeta)&=\frac{1-8(N-1)\zeta^2-\sqrt{1-16(N-1)\zeta^2}}{8(N-1)\zeta^2}
 =4(N-1)\zeta^2+O(|\zeta|^4) \label{eq:z5} \\
z_6=z_6(\zeta)&=\frac{1-8(N-1)\zeta^2+\sqrt{1-16(N-1)\zeta^2}}{8(N-1)\zeta^2}
 =\frac1{4(N-1)\zeta^2}+O(1) \label{eq:z6}
\end{align}
and the quantity~\eqref{eq:2ndterm} equals
\begin{equation}\label{eq:z3456int}
\frac1{8N\sqrt{N-1}\zeta^2\pi i}\int_{|z|=1}
 \frac{p(z)}{(z-z_3)(z_4-z)\sqrt{z(z-z_5)(z_6-z)}}\,dz.
\end{equation}
But we have
\[
\frac{p(z)}{(z-z_3)(z_4-z)}
 =4(N-1)\zeta^2+\frac{a_1}{z-z_3}+\frac{a_2}{z_4-z},
\]
where
\[
a_1=\frac{(N-2)^2\zeta^2z_3}{\sqrt{1-4N^2\zeta^2}},
\qquad a_2=\frac{(N-2)^2\zeta^2z_4}{\sqrt{1-4N^2\zeta^2}}.
\]
Thus,~\eqref{eq:z3456int} becomes
\begin{align}
&\frac{\sqrt{N-1}}{2N\pi i}\int_{|z|=1}\frac{dz}{\sqrt{z(z-z_5)(z_6-z)}} \label{eq:z56int} \\[2ex]
&+\frac{(N-2)^2 z_3}{8N\pi i\sqrt{N-1}\sqrt{1-4N^2\zeta^2}}
 \int_{|z|=1}\frac{dz}{(z-z_3)\sqrt{z(z-z_5)(z_6-z)}} \label{eq:z356int} \\[2ex]
&+\frac{(N-2)^2 z_4}{8N\pi i\sqrt{N-1}\sqrt{1-4N^2\zeta^2}}
 \int_{|z|=1}\frac{dz}{(z_4-z)\sqrt{z(z-z_5)(z_6-z)}}. \label{eq:z456int}
\end{align}
Using the contour in Figure~\ref{fig:contour}, but with $z_5$ replacing $z_1$,
we see that the term~\eqref{eq:z56int} equals
\begin{equation}\label{eq:z56intresult}
\frac{2\sqrt{N-1}}{N\pi\sqrt{z_6}}\int_0^{\pi/2}\frac{d\phi}{\sqrt{1-\frac{z_5}{z_6}\sin^2\phi}}
=\frac{2\sqrt{N-1}}{N\pi\sqrt{z_6(\zeta)}}F_1(\sqrt{\tfrac{z_5(\zeta)}{z_6(\zeta)}}),
\end{equation}
while
the term~\eqref{eq:z456int} equals
\begin{gather}
\notag
\frac{(N-2)^2}{2N\pi\sqrt{N-1}\sqrt{1-4N^2\zeta^2}\sqrt{z_6}}
\int_0^{\pi/2}\frac{d\phi}{(1-\frac{z_5}{z_4}\sin^2\phi)\sqrt{1-\frac{z_5}{z_6}\sin^2\phi}} \\
=\frac{(N-2)^2}{2N\pi\sqrt{N-1}\sqrt{1-4N^2\zeta^2}\sqrt{z_6(\zeta)}}
 \,\Pi_1(\tfrac{z_5(\zeta)}{z_4(\zeta)},\sqrt{\tfrac{z_5(\zeta)}{z_6(\zeta)}}),
\label{eq:z456intresult}
\end{gather}
where $\Pi_1$ is Legendre's complete elliptic integral of the third kind and where in~\eqref{eq:z56intresult}
and~\eqref{eq:z456intresult}, we take the branch of $\sqrt{z_6(\zeta)}$
that is close to $\frac1{2\sqrt{N-1}\zeta}$.
We now consider the term~\eqref{eq:z356int}.
We see from the asymptotics~\eqref{eq:z3}, \eqref{eq:z5} and~\eqref{eq:z6}
that for $\zeta$ sufficiently small we have $|z_5|<|z_3|<1$.
Hence, picking up the residue at $z_3$, we see that the term~\eqref{eq:z356int} equals
\begin{align}
&\frac{(N-2)^2z_3}{4N\sqrt{N-1}\sqrt{1-4N^2\zeta^2}\sqrt{z_3(z_3-z_5)(z_6-z_3)}} \label{eq:res} \\
\notag
&\qquad-\frac{(N-2)^2}{2N\pi\sqrt{N-1}\sqrt{1-4N^2\zeta^2}\sqrt{z_6}}
\int_0^{\pi/2}\frac{d\phi}{(1-\frac{z_5}{z_3}\sin^2\phi)\sqrt{1-\frac{z_5}{z_6}\sin^2\phi}} 
\end{align}
where we choose the branch of
\[
\sqrt{z_3(\zeta)(z_3(\zeta)-z_5(\zeta))(z_6(\zeta)-z_3(\zeta))}
=\sqrt{\tfrac{N^2(N-2)^2}{4(N-1)}\zeta^2+O(|\zeta|^4)}
\]
that is close to $\tfrac{N(N-2)}{2\sqrt{N-1}}\zeta$ and, again, the branch of $\sqrt{z_6(\zeta)}$
that is close to $\frac1{2\sqrt{N-1}\zeta}$.
However,
\[
(z_3-z_5)(z_6-z_3)=\frac{(N-2)^2}{4N^2(N-1)\zeta^2}z_3
\]
and the residue~\eqref{eq:res} equals
\[
\frac{(N-2)\zeta}{2\sqrt{1-4N^2\zeta^2}},
\]
which exactly cancels~\eqref{eq:2-N}.
Collecting all terms, we have
\begin{align*}
 C_{2,\ldots,2}(\zeta)=&\,\frac2{N\pi\sqrt{N-1}\sqrt{z_6(\zeta)}}\bigg(
(N-1)F_1(\sqrt{\tfrac{z_5(\zeta)}{z_6(\zeta)}}) \\
&\qquad+\frac{(N-2)^2}{4\sqrt{1-4N^2\zeta^2}}
 \big(\Pi_1(\tfrac{z_5(\zeta)}{z_4(\zeta)},\sqrt{\tfrac{z_5(\zeta)}{z_6(\zeta)}})
  -\Pi_1(\tfrac{z_5(\zeta)}{z_3(\zeta)},\sqrt{\tfrac{z_5(\zeta)}{z_6(\zeta)}})\big)\bigg),
\end{align*}
which yields for the Green function of the standard random walk on $G_{2,\ldots,2}$ the expression
\begin{align*}
\Gc_{2,\ldots,2}(z)=\frac8{N\pi\sqrt{w_6(z)}}\bigg(&
(N-1)F_1\bigg(\sqrt{\frac{w_5(z)}{w_6(z)}}\bigg) \\
&+\frac{(N-2)^2}{2\sqrt{4-N^2z^2}}
 \bigg(
\begin{aligned}[t]&\Pi_1\bigg(\frac{N^2w_5(z)}{(N-1)w_4(z)},
                   \sqrt{\frac{w_5(z)}{w_6(z)}}\bigg) \\
               &-\Pi_1\bigg(\frac{N^2w_5(z)}{(N-1)w_3(z)},
                   \sqrt{\frac{w_5(z)}{w_6(z)}}\bigg)\bigg)\bigg),
\end{aligned}
\end{align*}
where
\begin{align*}
w_3(z)&=8-N^2z^2-4\sqrt{4-N^2z^2} \\
w_4(z)&=8-N^2z^2+4\sqrt{4-N^2z^2} \\
w_5(z)&=2-(N-1)z^2-2\sqrt{1-(N-1)z^2} \\
w_6(z)&=2-(N-1)z^2+2\sqrt{1-(N-1)z^2}.
\end{align*}

\end{example}

\begin{example}\label{ex:23}
Consider the case $G=G_{2,3}$.
Using Maple to find a Groebner basis, one quickly computes $Q_{2,3}$ from $Q_2$ and $Q_3$.
(We checked this result in Mathematica
by back--substitution and elimination.)
We found that $Q_{2,3}=Q_{2,3}(R,b,\xi)$ from Proposition~\ref{prop:Z*ZRGr} is
\begin{align*}
Q_{2,3}&=b^5R^6 + 7b^4R^5 + b^3(19 - 12b^2 - 3b^2\xi)R^4 + 
 b^2(25 - 56b^2 - 14b^2\xi - 2b^3\xi)R^3 \\
 &\quad+b(16 - 93b^2 + 21b^4 - 23b^2\xi - 7b^3\xi + 12b^4\xi + 3b^4\xi^2 )R^2 \\
 &\quad+  (4 - 65b^2 + 49b^4 - 16b^2\xi - 9b^3\xi + 28b^4\xi - 6b^5\xi + 
   7b^4\xi^2 - 6b^5\xi^2)R \\
 &\quad - b(16 - 28b^2 + 2b^4 + 4\xi + 4b\xi - 17b^2\xi + 
   7b^3\xi - 3b^4\xi - 4b^2\xi^2 + 7b^3\xi^2  \\
 &\qquad\quad- b^4\xi^2 + b^4\xi^3).
\end{align*}
From this, we find $P_{2,3}=P_{2,3}( C,b,\xi)$ from Proposition~\ref{prop:Z*ZRGr} is
\begin{equation}\label{eq:P23}
\begin{aligned}
P_{2,3}&=(1 - 12b^2 + 21b^4 - 2b^6 - 3b^2\xi - 2b^3\xi + 12b^4\xi - 6b^5\xi + 3b^6\xi + 
   3b^4\xi^2 \\
&\qquad\quad - 6b^5\xi^2 + b^6\xi^2 - b^6\xi^3) C^3 \\
&\quad+ (b - 8b^3 + 7b^5 - 2b^3\xi - b^4\xi + 4b^5\xi - b^6\xi + b^5\xi^2 - b^6\xi^2) C^2 \\
&\quad- (b^2 - 3b^4 - b^4\xi + b^5\xi - b^6\xi) C
-b^3 + b^5.
\end{aligned}
\end{equation}
The $B$--valued Cauchy transform, $ C^{(B)}_{2,3}$, of the adjacency operator $T$
is an algebraic function of degree $3$.
To compute explicitly the $\Cpx$--valued Cauchy transform of $T$ by
integrating as was performed in Examples~\ref{ex:22} and~\ref{ex:2n2} seems, thus, to be difficult.
However, starting from $C^{(B)}_{2,3}(b)=b+O(\|b\|^2)$, further
terms of the power series expansion for $ C^{(B)}_{2,3}$ can be computed from the polynomial $P_{2,3}$.
We obtain, for instance,
\begin{equation}\label{eq:C23Bexp}
\begin{aligned}
 C^{(B)}_{2,3}(b)=&\,b+ (4 + \xi)b^3 + \xi b^4 + (26 + 12\xi + \xi^2)b^5 + 5\xi(3 + \xi) b^6  \\
& + ( 194 + 132\xi + 25\xi^2 + \xi^3) b^7+ 7\xi(5 + \xi)(5 + 2\xi) b^8 \\
& + (1542 + 1392\xi + 432\xi^2 + 52\xi^3 + \xi^4)b^9 \\
& + \xi(1887 + 1593\xi + 406\xi^2 + 30\xi^3) b^{10} \\
& + (12714 + 14320 \xi + 6275 \xi^2 + 1350 \xi^3 + 125 \xi^4 + \xi^5)b^{11} + O(\|b\|^{12}).
\end{aligned}
\end{equation}
Taking $\tau$ of~\eqref{eq:C23Bexp} and using
\[
\tau(\xi^n)=
\begin{cases}
\binom{n}{n/2}&n\text{ even} \\
0&n\text{ odd,}
\end{cases}
\]
we get the expansion for the $\Cpx$--valued Cauchy transform of $T$, which
gives the following expression for the first several
terms of the Green function for the standard random walk on $G_{2,3}$:
\begin{align*}
\Gc_{2,3}(z)=&\,1 + 4\,(\tfrac z4)^2 + 28\,(\tfrac z4)^4 + 10\,(\tfrac z4)^5 
 + 244\,(\tfrac z4)^6 + 210\,(\tfrac z4)^7 \\
& + 2412\,(\tfrac z4)^8 + 3366\,(\tfrac z4)^9 + 26014\,(\tfrac z4)^{10} + O(|z|^{11}).
\end{align*}
\end{example}

\begin{example}\label{ex:24}
Here is the case $G=G_{2,4}$, using notation as in the previous example:
\begin{align*}
P_{2,4}=&
\,(1 - 16 b^2 + 60 b^4 - 32 b^6 + 4 b^8 - 4 b^2 \xi + 
   30 b^4 \xi - 40 b^6 \xi + 4 b^8 \xi + 6 b^4 \xi^2 \\
&\quad - 28 b^6 \xi^2 - 3 b^8 \xi^2 - 4 b^6 \xi^3 - 
   2 b^8 \xi^3 + b^8 \xi^4) G^4 \\
&+ (2 b - 24 b^3 + 60 b^5 - 16 b^7 - 6 b^3 \xi + 30 b^5 \xi - 20 b^7 \xi + 6 b^5 \xi^2 \\
&\quad -   14 b^7 \xi^2 - 2 b^7 \xi^3) G^3 \\
& + (-2 b^4 + 10 b^6 + b^6 \xi - 2 b^8 \xi - 3 b^8 \xi^2) G^2  \\
& + (-2 b^3 + 8 b^5 - 2 b^7 + 2 b^5 \xi - b^7 \xi) G 
-b^4 + 2 b^6 \\
 C^{(B)}_{2,4}(b)=&\,
b + (4 + \xi)b^3  + (26 + 13 \xi + \xi^2)b^5  + (196 + 150 \xi + 30 \xi^2 + \xi^3)b^7 \\
&  + (1588 + 1644 \xi + 545 \xi^2 + 60 \xi^3 + \xi^4)b^9  \\
& +  (13424 + 17540 \xi + 8160 \xi^2 + 1585 \xi^3 + 110 \xi^4 + \xi^5)b^{11}+O(\|b\|^{13}) \\
\Gc_{2,4}(z)=&\,1 + 4 \,(\tfrac z4)^2 + 28 \,(\tfrac z4)^4 + 256 \,(\tfrac z4)^6 
 + 2684 \,(\tfrac z4)^8 + 30404 \,(\tfrac z4)^{10} + O(|z|^{12}).
\end{align*}
\end{example}

\begin{example}\label{ex:25}
Here is the case $G=G_{2,5}$.
\begin{align*}
P_{2,5}=&\,
(1 - 20 b^2 + 115 b^4 - 180 b^6 + 45 b^8 - 2 b^{10} - 
   5 b^2 \xi + 60 b^4 \xi - 2 b^5 \xi - 145 b^6 \xi  \\
&\quad - 30 b^7 \xi + 20 b^8 \xi + 10 b^9 \xi - 
   5 b^{10} \xi + 10 b^4 \xi^2 - 60 b^6 \xi^2 - 20 b^7 \xi^2 + 25 b^8 \xi^2 \\
&\quad - 10 b^9 \xi^2 + 
   b^{10} \xi^2 - 10 b^6 \xi^3 + 20 b^8 \xi^3 - 10 b^9 \xi^3 + 5 b^{10} \xi^3 + 5 b^8 \xi^4 - b^{10} \xi^5) G^5
 \displaybreak[0] \\
&+(3 b - 48 b^3 + 207 b^5 - 216 b^7 + 27 b^9 - 12 b^3 \xi + 
   108 b^5 \xi - 3 b^6 \xi - 174 b^7 \xi \\
&\quad - 27 b^8 \xi + 12 b^9 \xi + 3 b^{10} \xi + 18 b^5 \xi^2 - 
   72 b^7 \xi^2 - 18 b^8 \xi^2 + 15 b^9 \xi^2 - 3 b^{10} \xi^2 \\
&\quad - 12 b^7 \xi^3 + 12 b^9 \xi^3 - 
   3 b^{10} \xi^3 + 3 b^9 \xi^4) G^4  \displaybreak[0] \\
&+ (2 b^2 - 27 b^4 + 97 b^6 - 75 b^8 + 2 b^{10} - 6 b^4 \xi + 45 b^6 \xi - 6 b^7 \xi - 
   61 b^8 \xi \\
&\quad - 7 b^9 \xi - b^{10} \xi + 6 b^6 \xi^2 - 21 b^8 \xi^2 - 8 b^9 \xi^2 + 4 b^{10} \xi^2 - 
   2 b^8 \xi^3 + 3 b^{10} \xi^3) G^3  \displaybreak[0] \\
&+ (-2 b^3 + 13 b^5 - 7 b^7 - 5 b^9 + 4 b^5 \xi - 6 b^7 \xi - 4 b^8 \xi - 8 b^9 \xi \\
&\quad - 2 b^7 \xi^2 - 3 b^9 \xi^2 - b^{10} \xi^2) G^2  \displaybreak[0] \\
&+ (-3 b^4 + 15 b^6 - 11 b^8 + 3 b^6 \xi - 3 b^8 \xi - b^9 \xi - b^{10} \xi) G 
-b^5 + 3 b^7 - b^9  \displaybreak[1]\\
C^{(B)}_{2,5}(b)=&\,
b  + (4 + \xi) b^3 + (26 + 12 \xi + \xi^2) b^5 +  \xi b^6 + 
 (196 + 132 \xi + 24 \xi^2 + \xi^3) b^7 \\
& + (21 \xi + 7 \xi^2) b^8 + 
 (1590 + 1408 \xi + 400 \xi^2 + 40 \xi^3 + \xi^4) b^9 \\
& + (306 \xi + 189 \xi^2 + 27 \xi^3) b^{10} \\
& +  (13482 + 14800 \xi + 5741 \xi^2 + 940 \xi^3 + 60 \xi^4 + \xi^5) b^{11} \\
& + (3861 \xi + 3388 \xi^2 + 924 \xi^3 + 77 \xi^4)b^{12} + O(\|b\|^{13}) \\
\Gc_{2,5}(z)=&\,
1 + 4\,(\tfrac z4)^2 + 28\,(\tfrac z4)^4 + 244\,(\tfrac z4)^6 + 14\,(\tfrac z4)^7 
+ 2396\,(\tfrac z4)^8  \\
&+ 378\,(\tfrac z4)^9 + 25324\,(\tfrac z4)^{10} 
+ 7238\,(\tfrac z4)^{11}+O(|z|^{12}).
\end{align*}
\end{example}

\bibliographystyle{plain}

\end{document}